\documentclass{amsart}
\usepackage{amssymb,amscd,enumerate}
\usepackage[french]{babel}
\usepackage[latin1]{inputenc}
\usepackage[T1]{fontenc}
\usepackage[11pt,curve,matrix,arrow,frame,graph]{xy}
\usepackage{slashed}
\newcounter{intro}

\newtheorem{theo}[intro]{Th\'eor\`eme}

\newtheorem{thm}{Th\'eor\`eme}[section]

\newtheorem{prop}[thm]{Proposition}

\theoremstyle{remark}

\newtheorem*{merci}{Remerciements}

\numberwithin{equation}{section}
\newcounter{counteroman}

\newcommand{\cref}[1]{Corollary~\ref{#1}}

\newcommand{\R}{\mathbb{R}}

\newcommand{\bS}{\mathbb{S}}

\newcommand{\mcO}{\mathcal{O}}

\let\ve=\varepsilon

\DeclareMathOperator{\vol}{vol}

\DeclareMathOperator{\scal}{Scal}
\DeclareMathOperator{\ric}{Ricci}

\DeclareMathOperator{\eucl}{eucl}

\begin{document}

\title[Inégalité de Sobolev et volume asymptotique]
{Inégalité de Sobolev et volume asymptotique}
\author{Gilles Carron}
\address{Laboratoire de Math\'ematiques Jean Leray (UMR 6629), 
Universit\'e de Nantes,
2, rue de la Houssini\`ere, B.P.~92208, 44322 Nantes Cedex~3, France}
\email{Gilles.Carron@math.univ-nantes.fr}


\date{\today}
\begin{abstract}En 1999, M. Ledoux a démontré qu'une variété complète à courbure de Ricci positive ou nulle qui vérifiait une inégalité de Sobolev euclidienne était euclidienne. On présente un raccourci de la preuve. De plus ces arguments permettent  un raffinement d'un résultat de B-L. Chen et X-P. Zhu à propos des variétés localement conformément plate à courbure de Ricci positive ou nulle. Enfin, on étudie ce qui se passe lorsque l'hypothèse sur  la courbure de Ricci est remplacée par une hypothèse sur la courbure scalaire.\end{abstract}
\maketitle
\section{Introduction}
L'objectif de cette note est d'abord de pr\'esenter une l\'eg\`ere simplification de la preuve d'un tr\`es joli r\'esultat de M. Ledoux (\cite{Ledoux}). Pour pr\'esenter ce r\'esultat, nous introduisons les meilleurs constantes $K(n,p)$ ($p\in [1,n[$) de l'in\'egalit\'e de Sobolev euclidienne : 
\begin{equation}\label{Sobeu}\forall f\in C^\infty_0(\R^n)\, ,\, \|f\|^p_{\frac{pn}{n-p}}\le K(n,p)^p \,\|df\|_p^p\,\,;\end{equation}
c'est à dire
$$K(n,p)^{-1}=\inf_{ f\in C^\infty_0(\R^n)}\frac{\|df\|_p}{\,\,\,\,\, \|\,f\,\|_{\frac{pn}{n-p}}}$$
Grâce aux travaux de G.Talenti et T. Aubin (\cite{Aubin}, \cite{Talenti}), on connaît la valeur de $K(n,p)$ et on sait de plus, que pour $p\in ]1,n[$, les fonctions 
$$f_\lambda(x)=\frac{1}{\left(\lambda^q+\|x\|^q\right)^{\frac{n-p}{p}}} ,\,\, \mathrm{avec}\,\, q=\frac{p}{p-1}\,\,$$
réalisent l'égalité dans l'inégalité de Sobolev euclidienne (\ref{Sobeu}). Le résultat de M. Ledoux est le suivant
\begin{thm}
Soit $(M^n,g)$ une variété riemannienne complète à courbure de Ricci positive ou nulle. Si pour un $p\in [1,n[$ celle ci vérifie l'inégalité de Sobolev :
$$\forall f\in C^\infty_0(M)\,,\,\, \|f\|^p_{\frac{pn}{n-p}}\le K(n,p)^p\,\|df\|_p^p\,\,$$
alors 
$(M^n,g)$ est isométrique à l'espace euclidien $\R^n$.
\end{thm}

Un raisonnement basé sur le caractère infinitésimallement euclidien de la géométrie  riemannienne montre facilement que si 
une variété riemannienne $(M^n,g)$ vérifie l'inégalité de Sobolev 
\begin{equation}\label{SobM}\forall f\in C^\infty_0(M)\,,\,\, \|f\|^p_{\frac{pn}{n-p}}\le A^p \,\|df\|_p^p\,\,\end{equation}
alors $$A\ge K(n,p)\,\,.$$

On sait de plus qu'une telle inégalité de Sobolev (\ref{SobM}) implique en toute généralité une croissance euclidienne pour le volume des boules géodésiques (\cite{Aka}, \cite{C} ,\cite[theorem 3.1.5]{SC}) :
$$\forall x\in M, \forall r>0\,\,,\,\,\vol B(x,r)\ge C(p,n) \left(\frac{r}{A}\right)^n\,\,.$$

Lorsque la courbure de Ricci de la métrique $g$ est positive ou nulle, le théorème de comparaison de Bishop-Gromov nous apprend que le quotient $$\frac{\vol B(o,r)}{r^n}$$ est une fonction décroissante de $r$, on peut alors définir
$$\lim_{r\to +\infty} \frac{\vol B(o,r)}{r^n}=:\nu_{\infty}$$ le volume asymptotique de $(M^n,g)$. On a donc
$\nu_{\infty}\le \omega_n$
où $\omega_n$ est le volume de la boule euclidienne de rayon $1$. De plus l'égalité 
$\nu_{\infty}= \omega_n$ implique que la variété $(M^n,g)$ est isométrique à l'espace euclidien $\R^n$.

Le but de la preuve de M. Ledoux consiste à démontrer que l'hypothèse de positivité de la courbure de Ricci  et l'inégalité de Sobolev 
$$\forall f\in C^\infty_0(M)\,\,\, \|f\|^p_{\frac{pn}{n-p}}\le K(n,p)^p\,\|df\|_p^p\,\,$$
implique que $\nu_{\infty}\ge \omega_n\,\,.$
Ceci se fait par l'étude d'une inéquation différentielle impliquée par l'inégalité de Sobolev.

La preuve fournie par M. Ledoux a été reprise par de nombreux auteurs et pour différentes inégalités de type Sobolev pour démontrer que si la constante $A$ de l'inégalité de Sobolev (\ref{SobM}) est presque euclidienne et si la courbure de Ricci est positive ou nulle alors
le volume asymptotique est presque maximal et donc d'après un résultat de T. Colding la variété est difféomorphe à $\R^n$
( \cite{LX}, \cite{CX}, \cite{RC}, \cite{LX}, \cite{X1}, \cite{X2}, \cite{X3}). Dans un papier récent S. Pigola et  G. Veronelli (\cite{PV})  utilise une autre méthode basée sur des théorèmes de comparaison pour le laplacien.

Ici, on va démontrer le résultat suivant 
\begin{theo}\label{thm}Soit $(M^n,g)$ une variété Riemannienne complète qui vérifie l'inégalité de Sobolev :
$$\forall f\in C^\infty_0(M)\,,\,\, \|f\|^p_{\frac{pn}{n-p}}\le A^p \,\|df\|_p^p$$ 
si le quotient $$\frac{\vol B(o,r)}{r^n}$$ admet une limite lorsque $r$ tend vers l'infini :
$$\nu_\infty=\lim_{r\to +\infty} \frac{\vol B(o,r)}{r^n}$$ alors
$$\nu_\infty\ge \omega_n \left(\frac{K(n,p)}{A}\right)^{\frac{1}{n}}\,\,.$$
\end{theo}

Remarquons que la preuve de M. Ledoux pourrait être adaptée pour démontrer ce résultat ;  cf. la preuve de \cite[theorem 3.3]{LX}. De plus, il sera évident que ce résultat et sa preuve peuvent être généralisé à beaucoup d'autres inégalités de type Sobolev (voir par exemple \cite{CX},\cite{X2},\cite{X3}). 

Pour les variétés à courbure de Ricci presque positive, on sait que ce quotient a une limite. C'est à dire on suppose que pour un point $o\in M$ on a la minoration :
$$\ric_g\ge -(n-1) G(d(o,\,.\,)) g$$
où $G$ est une fonction continue positive sur $[0,+\infty[$ qui vérifie 
$$b:=\int_0^{+\infty}G(t)dt<+\infty\,\,.$$
Les théorèmes de comparaison montrent que si $h$ est la fonction solution du problème de Cauchy
$$h''=Gh,\, h(0)=0, h'(0)=1$$
alors 
pour $$V(r)=\int_0^r h(t)^{n-1} dt\,\,,$$
le quotient $$\frac{\vol B(o,r)}{V(r)}$$ est une fonction décroissante.
De plus  la fonction $h$ est convexe, le quotient $h(r)/r$ est une fonction croissante et l'on a
pour $$\alpha:=\lim_{r\to\infty} \frac{h(r)}{r}\,\,\, \mathrm{alors}\,\,\,1+b\le \alpha\le e^b\,\,.$$
Ce qui montre que 
dans ce cas le quotient 
$\frac{\vol B(o,r)}{r^n}$ a bien une limite en $+\infty$ égale à
$$\lim_{r\to +\infty} \frac{\vol B(o,r)}{r^n}=\frac{\,\,\,\alpha^{n-1}}{n}\,\, \lim_{r\to\infty}\frac{\vol B(o,r)}{V(r)}\,\,.$$
Ainsi, notre théorème \ref{thm} permet de retrouver des résultats récents de S. Pigola et  G. Veronelli (\cite{PV}) et de L. Adriano et C. Xia (\cite{LX}). 

Par ailleurs, notre preuve nous permettra également de re-démontrer un résultat de rigidité de B-L. Chen et X-P. Zhu (\cite{ChenZhu})
\begin{thm} Soit $(M^n,g)$ une variété riemannienne complète non compacte localement conformément plate à courbure de Ricci positive ou nulle si
$$\lim_{r\to +\infty} \frac{r^2}{\vol B(o,r)} \int_{B(o,r)} \scal_gd\!\vol_g\,=0$$ alors
$(M^n,g)$ est isométrique à l'espace euclidien $\R^n$.
\end{thm}
En fait nous démontrerons un petit raffinement de ce résultat à savoir :
\begin{theo} Soit $(M^n,g)$ une variété riemannienne complète non compacte localement conformément plate à courbure de Ricci positive ou nulle alors soit la  variété est plate soit nous avons l'inégalité :
$$\omega_n^{\frac{2}{n}}-\nu_\infty^{\frac{2}{n}}\le n^{-2}\,\,\limsup_{r\to +\infty} \frac{1}{\left(\vol B(o,r)\right)^{\frac{n-2}{n}}} \int_{B(o,r)} \scal_gd\!\vol_g\,\,\,.$$ 
\end{theo}

Enfin, on remarque qu'une variété riemannienne $(M^n,g)$ complète à courbure scalaire positive ou nulle vérifiant l'inégalité de Sobolev optimale :
$$\forall f\in C^\infty_0(M)\,, \,\, \big\|f\big\|^2_{\frac{2n}{n-2}}\le K(n,2)^2 \,\big\|df\big\|_2^2\,\,,$$ vérifie aussi l'inégalité suivante
$$\forall f\in C^\infty_0(M)\,, \,\, \big\|f\big\|^2_{\frac{2n}{n-2}}\le K(n,2)^2 \,\int_M \left[ |df|^2+\frac{n-2}{4(n-1)} \scal_g f^2\right]d\!\vol_g\,\,\,.$$
Ainsi l'invariant de Yamabe de la métrique $g$ est égale à celui de la sphère ronde :
$$Y([g]):=\inf_{f\in C^\infty_0(M)}\frac{\int_M \left[ |df|^2+\frac{n-2}{4(n-1)} \scal_g f^2\right]d\!\vol_g}{\big\|f\big\|^2_{\frac{2n}{n-2}}}=Y(\bS^n).$$

Grâce au travaux de T.Aubin et R. Schoen (\cite{Au2}, \cite{Schoen}), on sait que pour une variété compacte $(M^n,g)$ non conformément équivalente à la sphère ronde
on a 
$$Y([g])<Y(\bS^n)\,\,.$$
Ce résultat est fondamental pour trouver une métrique conforme à courbure scalaire constante et il permet de conclure le programme initié
par H.Yamabe et N. Trudinger (\cite{Yam},\cite{Tr}).
Ainsi une généralisation du résultat de M. Ledoux serait de déterminer les variétés riemanniennes complètes à courbure scalaire positive ou nulle dont l'invariant de Yamabe est égale à celui de la sphère ronde. On remarque qu'il y a  de nombreux exemples de telles variétés, en effet selon \cite[Prop 2.2]{SY}, nous avons 
\begin{prop} Soit $(M^n,g)$ une variété riemannienne  complète telle qu'il existe une immersion conforme 
$$\Phi\colon (M,g)\rightarrow (\bS^n,can)$$ alors
$$Y([g])=Y(\bS^n)\,\,.$$
\end{prop}
Le résultat spectaculaire de l'article de R.Schoen et S-T.Yau est qui si de plus la courbure scalaire est positive ou nulle alors 
(dans la plupart des cas) une telle immersion est injective.
Notre résultat est le suivant :
\begin{theo}Soit $(M,g)$ une variété riemannienne complète non compacte à courbure scalaire positive :
$$\scal_g\ge 0$$ et dont l'invariant de Yamabe est celui de la sphère ronde 
$$Y([g])=Y(\bS^n)\,\,.$$
Supposons l'une des hypothèses suivantes vérifiées :
\begin{enumerate}[i)]
\item $n=\dim M\ge 6$,
\item $M$ est spin \footnote{Par exemple si $M$ est orienté et si $n=3$.},
\item $M^n$ est localement conformément plate de dimension $n\not=5$.
\end{enumerate}
alors $(M,g)$ est conformément équivalente à un ouvert $\Omega\subset \bS^n$ telle que
$$\dim_H(\bS^n\setminus \Omega)\le (n-2)/2.$$
\end{theo}
La preuve de ce théorème est basée sur l'analyse de Schoen-Yau et sur une adaptation pour les variétés complètes de la preuve du résultat de Aubin et Schoen, en particulier on utilisera une nouvelle fonction test qui est inspirée de l'analyse faite par R. Schoen.
Il faut ici remarquer que si le théorème de la masse positive généralisée est vraie (cf. le discours de R. Schoen et S-T. Yau pour énoncer \cite[prop 4.4']{SY}) alors la conclusion est toujours vraie sans les hypothèses i) ou ii) ou iii).

\begin{merci}
Je remercie  V. Minerbe  pour ses commentaires avisés et je suis aussi partiellement financé par le projet ACG: ANR-10-BLAN 0105
\end{merci}
\section{Preuve du théorème A :}
Le cas $p=1$ est relativement aisé car l'inégalité de Sobolev 
$$\forall f\in C^\infty_0(M)\,\,\, \|f\|_{\frac{n}{n-1}}\le A \,\|df\|_1$$ est équivalente à l'inégalité isopérimétrique 
$$\forall \Omega\subset M,\,\, A\vol(\partial\Omega)\ge \vol (\Omega)^{\frac{n-1}{n}}$$
qui, elle, implique la minoration 
$$\vol B(o,r)\ge  \left(\frac{r}{n\, A}\right)^n\,\,.$$

Concernant le cas où $p\in ]1,n[$, on introduit, comme M. Ledoux, les fonctions 
$$f_\lambda(x)=\frac{\lambda^{\frac{(q-1)(n-p)}{p}}}{ \left(\lambda^q+d(o,x)^q\right)^{\frac{n-p}{p}}}.$$
où $q=p/(p-1)$. On notera 
$$\theta(r)=\frac{\vol B(o,r)}{r^n}\,\,.$$
En reprenant les calculs fait par M. Ledoux, on a facilement :
$$ \|f_\lambda\|_{\frac{pn}{n-p}}^p=\left(\int_0^{+\infty} \frac{nq r^{q-1+n}\lambda^{n(q-1)}}{\left(\lambda^q+r^q\right)^{n+1}}  \theta(r) dr\right)^{1-\frac{p}{n}}$$
En faisant alors le changement de variables $r=\lambda\rho$ on obtient :
$$ \|f_\lambda\|_{\frac{pn}{n-p}}^p=\left(\int_0^{+\infty} \frac{nq \rho^{q-1+n}}{\left(1+\rho^q\right)^{n+1}}  \theta(\lambda\rho) d\rho\right)^{1-\frac{p}{n}}$$
En notant 
$$A_{n,p}=\left(\int_0^{+\infty} \frac{nq \rho^{q-1+n}}{\left(1+\rho^q\right)^{n+1}}  d\rho\right)^{1-\frac{p}{n}}$$  
grâce au théorème de convergence dominée, on obtient donc :
\begin{equation}
\label{equa1}
\lim_{\lambda\to +\infty} \|f_\lambda\|_{\frac{pn}{n-p}}^p=  \left(\nu_\infty\right)^{1-\frac{p}{n}} A_{n,p}\,\,.
\end{equation}
De même, en introduisant la mesure de Stieljes : $dv(r)$ où $v(r)=\vol B(o,r)$, on trouve 
\begin{equation*}\begin{split}
 \|df_\lambda\|_{p}^p&=q^p\left(\frac{n-p}{p}\right)^p\,\int_M \frac{\lambda^{(q-1)(n-p)}d(o,x)^{(q-1)p}}{ \left(\lambda^q+d(o,x)^q\right)^{n}} d\!\vol_g(x)\,\\
 &=q^p\left(\frac{n-p}{p}\right)^p\,\int_M \frac{\lambda^{(q-1)(n-p)}r^{(q-1)p}}{ \left(\lambda^q+r^q\right)^{n}} dv(r)\,\\
 &=q^p\left(\frac{n-p}{p}\right)^p\lambda^{(q-1)(n-p)}\int_0^{+\infty} qr^{q-1}\left( \frac{nr^q}{\left(\lambda^q+r^q\right)^{n+1}}-\frac{1}{\left(\lambda^q+r^q\right)^n}\right) v(r) dr\\
 &=q^p\left(\frac{n-p}{p}\right)^p\int_0^{+\infty} q\rho^{n+q-1}\left(\frac{n\rho^q}{\left(1+\rho^q\right)^{n+1}}-\frac{1}{\left(1+\rho^q\right)^n}\right) \theta(\lambda\rho)d\rho
 \end{split}\end{equation*}
 D'où  en notant 
 \begin{equation*}\begin{split}
 B_{n,p}&=  \,\,q^p\left(\frac{n-p}{p}\right)^p\int_0^{+\infty} q\rho^{n+q-1}\left(\frac{n\rho^q}{\left(1+\rho^q\right)^{n+1}}-\frac{1}{\left(1+\rho^q\right)^n}\right) d\rho\,\,,\\
 &=  \,\,q^p\left(\frac{n-p}{p}\right)^p\, n\,\int_0^{+\infty}\frac{\rho^{n+q-1}}{\left(1+\rho^q\right)^{n}}d\rho\,\,,
 \end{split}\end{equation*}
 on obtient :
 \begin{equation}
\label{equa2}
\lim_{\lambda\to +\infty} \|df_\lambda\|_{p}^p=  \nu_\infty B_{n,p}\,\,. \end{equation}
Ainsi l'inégalité de Sobolev: 
$$\forall f\in C^\infty_0(M)\,\,\, \|f\|^p_{\frac{pn}{n-p}}\le A^p \,\|df\|_p^p$$
implique avec les égalités (\ref{equa1}) et (\ref{equa2}) :
$$\left(\nu_\infty\right)^{1-\frac{p}{n}} A_{n,p}\le  A^p  \nu_\infty B_{n,p} $$
Et puisque lorsqu'on effectue les calculs sur l'espace euclidien $\R^n$ on a égalité, on en déduit
$$\left(\nu_\infty\right)^{-\frac{p}{n}}\le A^p \frac{B_{n,p} }{A_{n,p}}=\left(\omega_n\right)^{-\frac{p}{n}}\left(\frac{A}{K(n,p)}\right)^p.$$
\section{Preuve du théorème B :}
On va maintenant démontrer le théorème B. Pour cela on se sert de la classification des variétés localement conformément plates à courbure de Ricci positive ou nulle obtenue dans (\cite{CH}). Cette classification implique qu'une telle variété, si elle n'est pas compacte, est 
\begin{itemize}
\item soit plate, 
\item soit isométrique à un quotient du produit $\R\times \bS^{n-1}$
\item soit globalement conformément équivalente à $\R^n$.\end{itemize}
Dans le deuxième cas, le volume asymptotique est nulle et la limite du terme de gauche est l'infini donc l'inégalité est bien valide.
Il faut donc traiter le cas où $M=\R^n$ équipé d'une métrique
$$g=u^{\frac{4}{n-2}}\,\mathrm{eucl}\,\,.$$
Dans ce cas, l'invariance conforme de l'invariant de Yamabe fournit l'inégalité de type Sobolev :
$$\forall f\in C^\infty_0(\R^n)\,\,\, \big\|f\big\|^2_{\frac{2n}{n-2}}\le K(n,2)^2 \,\left[\big\|df\big\|_2^2+\frac{n-2}{4(n-1)}\int_{\R^n} \scal_g f^2d\!\vol_g\,\right]\,\,.$$
On utilise les mêmes arguments en testant cette inégalité à la fonction
$$f_\lambda(x)=\frac{\lambda^{\frac{n-2}{2}}}{\left(\lambda^2+d(o,x)^2\right)^{\frac{n-2}{2}}}\,$$
pour obtenir
$$\left(\nu_\infty\right)^{1-\frac{2}{n}} A_{n,2}\le  K(n,2)^2\left(  \nu_\infty B_{n,2} +I\right)\,,$$
 où
$$I=\limsup_{\lambda\to \infty} \frac{n-2}{4(n-1)}\int_{\R^n} \scal_g f_\lambda^2\,d\!\vol_g\,.$$
Or, en introduisant la mesure de Stieljes $dS(r)$ où $S(r)=\int_{B(o,r)} \scal_gd\!\vol_g\,$,  on a 
\begin{equation*}\begin{split}
\int_{\R^n} \scal_g f_\lambda^2(x)\vol_g\,&=\int_{\R^n} \frac{\lambda^{n-2}}{\left(\lambda^2+d(o,x)^2\right)^{n-2}}\scal_g(x)d\!\vol_g(x)\,\\
&=\int_0^{+\infty} \frac{\lambda^{n-2}}{\left(\lambda^2+r^2\right)^{n-2}}dS(r)\\
&=\int_0^{+\infty} 2(n-2)\frac{\lambda^{n-2}r}{\left(\lambda^2+r^2\right)^{n-1}}S(r)dr
 \end{split}\end{equation*} 
 D'où en posant 
 $$\Sigma(r)=r^{2-n}\int_{B(o,r)} \scal_gd\!\vol_g\,\,\,,$$
\begin{equation*}\begin{split} \int_{\R^n} \scal_g f_\lambda^2d\!\vol_g\,&=2(n-2)\int_0^{+\infty} \frac{\lambda^{n-2}r^{n-1}}{\left(\lambda^2+r^2\right)^{n-1}}\Sigma(r)dr\\
&=2(n-2)\int_0^{+\infty} \frac{\rho^{n-1}}{\left(1+\rho^2\right)^{n-1}} \Sigma(\lambda\rho)d\rho\,\,
  \end{split}\end{equation*}
  D'où 
  $$I\le C_n \limsup_{r\to +\infty} r^{2-n}\int_{B(o,r)} \scal_gd\!\vol_g\,$$
  où
  $$C_n=\frac{(n-2)^2}{2(n-1)}\int_0^{+\infty} \frac{\rho^{n-1}}{\left(1+\rho^2\right)^{n-1}}d\rho\,\,.$$
  Et on obtient finalement :
$$\left(\nu_\infty\right)^{1-\frac{2}{n}} \frac{A_{n,2}}{B_{n,2}}\le  K(n,2)^2\left(  \nu_\infty + \frac{C_n}{B_{n,2}} \limsup_{r\to +\infty} r^{2-n}\int_{B(o,r)} \scal_gd\!\vol_g\,\right)$$
soit encore
$$\omega_n^{\frac{2}{n}}\le \nu_\infty^{\frac{2}{n}}+\frac{C_n}{B_{n,2}}\limsup_{r\to +\infty}\frac{\int_{B(o,r)} \scal_gd\!\vol_g\,}{\left(\nu_\infty r^n\right)^{1-\frac{2}{n}}}\,\,.$$
Pour finir on compare les deux constantes $C_n$ et $B_{n,2}$.
Rappelons que 
\begin{equation*}
\begin{split}
B_{n,2}&=(n-2)^2n\int_0^\infty\frac{\rho^{n+1}}{(1+\rho^2)^n}d\rho\\
&=-(n-2)^2n\int_0^\infty\left(\frac{d}{d\rho}\frac{1}{(1+\rho^2)^{n-1}}\right)\, \frac{\rho^n}{2(n-1)}d\rho\\
&=(n-2)^2n\int_0^\infty\frac{1}{(1+\rho^2)^{n-1}}\,\frac{n\rho^{n-1}}{2(n-1)}d\rho\\
&=n^2 C_n.
\end{split}\end{equation*}

\section{Preuve du théorème C}
La preuve de ce résultat se déroule avec les mêmes arguments que la preuve du résultat sus mentionné de T.Aubin et R.Schoen avec quelques modifications dues à la non compacité de la variété. On veut implémenter dans l'inégalité de Sobolev-Yamabe :
\begin{equation}\label{YS}\begin{split}
\, &\forall f\in C^\infty_0(M)\,\,\, \\
&\left(\int_M |f|^{\frac{2n}{n-2}} d\!\vol_g\right)^{1-\frac2n}\le  K(n,2)^2 \,\int_M \left[ |df|^2+\frac{n-2}{4(n-1)} \scal_g f^2\right]d\!\vol_g
\end{split}
\end{equation}
une bonne fonction test qui nous permettra de se ramener au résultat de Schoen-Yau.

\subsection{Estimées préliminaires}

Le point de départ est que cette inégalité de Sobolev (\ref{YS}) assure l'existence d'un noyau de Green minimal positif , noté $G(x,y)=G_x(y)$, pour l'opérateur de Yamabe :
$$L:=\Delta_g+\frac{n-2}{4(n-1)} \scal_g\,\,.$$
Pour alléger les notations on notera $$a_n= \frac{n-2}{4(n-1)}.$$
De plus si $x$ est un point de $M$ et $\mcO$ un voisinage compact de $x$ alors on sait que (cf. \cite[Cor. 2.3]{SY})
\begin{equation}\label{fini}\begin{split}
&\int_{M\setminus \mcO} G_x^{\frac{2n}{n-2}}d\!\vol_g<+\infty\\
 \mathrm{et }\,\,\,& 
\int_{M\setminus \mcO} \left[|dG_x|^2+a_n \scal_g G_x^2\right]\,d\!\vol_g<+\infty
\end{split}\end{equation}

On sait aussi que l'inégalité de Sobolev (\ref{YS}) implique une estimée gaussienne sur le noyau de la chaleur de $L$ : il y a une constante $C_n$ telle que 
$$\forall x,y\in M, \forall t>0 \,\,:\,\, e^{-tL}(x,y)\le C_n\, \frac{e^{-\frac{d^2(x,y)}{5t}}}{t^{\frac n2}},$$
ainsi on a 
$$G_x(y)=\int_0^{+\infty} e^{-tL}(x,y)dt\le \frac{C'_n}{d(x,y)^{n-2}}.$$

Notons $\Omega:=\{y\in M, G_x(y)< t\}$, c'est le complémentaire d'un voisinage compact de $x$ et en intégrant par parties, nous obtenons :
\begin{equation*}\begin{split}
\int_{\Omega_t}\left[ |dG_x|^2+a_n \scal_g G_x^2\right]d\!\vol_g&=\int_{\{G_x=t\}} G_x\frac{\partial G_x}{\partial n} \\
&=t\int_{\{ G_x=t\}} \frac{\partial G_x}{\partial n}\\
&=t\int_{\{G_x=t\}} \Delta G_x\\
&=t-t\int_{\{G_x>t\}} a_n \scal_g G_x\,d\!\vol_g
\end{split}
\end{equation*}
où $n$ désigne la normale unitaire sortante à $\partial \Omega_t$ et l'expression 
$\int_{\{G_x>t\}} \Delta G_x$ est à comprendre au sens des distributions.
On applique alors l'inégalité de Sobolev (\ref{YS}) à la fonction 
$\min(G_x,t)$ pour obtenir 
$$t^2\vol\{G_x>t\}^{1-\frac{2}{n}}\le K(n,2)^2\left[ t+t\int_{\{G_x>t\}} a_n \scal_g (t-G_x)d\!\vol_g
\right]\le K(n,2)^2\, t\,.$$
Ce qui permet d'obtenir l'estimation suivante 
\begin{equation}\label{volgro}
\vol\{G_x>t\}\le C t^{-\frac{n}{n-2}}.
\end{equation}

On remarque maintenant que l'invariance conforme du noyau de Green fait que ces estimées sont encore valides
pour une déformation conforme à support compact de $g$.
\subsection{Une bonne fonction test}
\subsubsection{Implémentation}
On fixe maintenant $x\in M$ et suivant \cite[thm. 5.1]{LP}, on modifie conformément $g$ dans un voisinage compact de $x$
pour obtenir une métrique $$\bar g=e^v g$$
où $ v(x)=1$, $v$ est à support dans un voisinage compact de $x$ et  
en $x$ toutes les dérivées covariantes symétriques du tenseur de Ricci  de $\bar g$ sont nulles jusqu'à l'ordre $N$ ($N$ étant choisi assez grand). Ceci implique qu'en coordonnées normales autour de $x$ on a :
$$d\!\vol_{\bar g}=\left(1+O\left(s^N\right)\right)s^{n-1}dsd\sigma $$
où $d\sigma$ est la mesure de Lebesgue de la sphère unité dans l'espace tangent $(T_xM, \bar g_x=g_x)$.

Dans ce cas on sait aussi que 
\begin{equation}\label{Deltascal}
\Delta_{\bar g} \scal_{\bar g}(x)=\frac{|W|^2(x)}{6}.
\end{equation}
où $|W|(x)$ est la norme du tenseur de Weyl de $g$ en $x$ (c'est aussi la norme du tenseur de Weyl de $\bar g$ en $x$).
Tous les calculs de cette sous-section (4.2) seront fait pour la métrique $\bar g$.

Posons $\sigma_{n-1}=n\omega_n$ le volume de la sphère unité de dimension $n-1$ et 
$$\alpha_n=\frac{1}{\sigma_{n-1}(n-2)}.$$
Soit $\lambda>0$, la fonction test que nous choisissons est 
$$f_\lambda:= \frac{\lambda^{\frac{n-2}{2}}}{\left(\lambda^2+r^2\right)^{\frac{n-2}{2}}}$$
où la fonction $r\colon M\rightarrow \R_+$ est définie par
$$G_x(y)=\frac{\alpha_n}{r(y)^{n-2}},$$
où on le rappelle $G_x$ est la fonction de Green conforme pour la métrique $\bar g$ avec pôle en $x$.
L'estimée (\ref{volgro}) assure que le volume des sous-lignes de niveau de $r$ croit au maximum de façon euclidienne :
$$\vol\{r\le R\}\le C R^n.$$
 Ainsi on a $f_\lambda\in L^{\frac{2n}{n-2}}$ et le même calcul que précédemment montre que 
\begin{equation}\begin{split}
& \|f_\lambda\|_{\frac{2n}{n-2}}^{\frac{2n}{n-2}}=\int_0^{+\infty}
  \frac{2n \rho^{1+n}}{\left(1+\rho^2\right)^{n+1}}  \theta(\lambda\rho) d\rho\\
  & \mathrm{avec} \,\,\theta(\tau)=\frac{\vol\{r\le \tau\}}{\tau^n}
  \end{split}
  \end{equation}
  
Ensuite on calcule comme précédemment :
\begin{equation*}\begin{split}
\int_M |df_\lambda|^2&=\lambda^{n-2}(n-2)^2\int_M \frac{r^2}{\left(\lambda^2+r^2\right)^n} |dr|^2 d\!\vol_{\bar g}\\
&=\lambda^{n-2}(n-2)^2\int_0^{+\infty} \frac{\tau^2}{\left(\lambda^2+\tau^2\right)^n}\left(\int_{r=\tau}  |dr|\right)d\tau\end{split}\end{equation*}
 Cependant 
 \begin{equation*}\begin{split}
 \int_{\{r=\tau\}}  |dr|&=\frac{\tau}{n-2} \int_{\{G_x=\alpha_n/\tau^n\}} \frac{|dG_x|}{G_x}\\
 &=\frac{\tau^{n-1}}{(n-2)\alpha_n}  \int_{\{G_x=\alpha_n/\tau^n\}} |dG_x|\\
 &=\frac{\tau^{n-1}}{(n-2)\alpha_n}\left[1- a_n\int_{\{r<\tau\}} \scal_g G_xd\!\vol_g\right]
 \end{split}\end{equation*}
 D'où en posant 
 $$v(\tau)=a_n\int_{\{r<\tau\}} \scal_g G_xd\!\vol_g$$ on obtient finalement :
 $$\int_M |df_\lambda|^2=\sigma_{n-1}\lambda^{n-2}(n-2)^2\left[\int_0^{+\infty}\frac{\tau^{n+1}}{\left(\lambda^2+\tau^2\right)^n}d\tau
  -\int\frac{\tau^{n+1}}{\left(\lambda^2+\tau^2\right)^n} v(\tau) d\tau\right].$$

On a aussi
\begin{equation*}\begin{split}
\int_{M} a_n \scal_g f_\lambda^2 d\!\vol_{\bar g}&= \int_M \frac{(\lambda r)^{n-2}}{\alpha_n \left(\lambda^2+r^2\right)^{n-2}}a_n\scal_{\bar g}G_xd\!\vol_{\bar g}\\
&=\int_0^{+\infty} \frac{(\lambda \tau)^{n-2}}{\alpha_n \left(\lambda^2+\tau^2\right)^{n-2}}dv(\tau)\\
&=\lambda^{n-2}(n-2)^2\sigma_{n-1} \int_0^{+\infty} \frac{\tau^4-\lambda^4}{\left(\lambda^2+\tau^2\right)^3} \left(\frac{\tau}{\lambda^2+\tau^2}\right)^{n-3}v(\tau)d\tau
 \end{split}\end{equation*}
 où comme précédemment nous avons utilisé la mesure de Stieljes $dv(\tau)$.
Au final il nous reste :
$$\int_M\left[ |df_\lambda|^2+ a_n \scal_g f_\lambda^2 \right]d\!\vol_{\bar g}=\omega_n B_{n,2}-J(\lambda)$$
où
$$\omega_nB_{n,2}=\sigma_{n-1}(n-2)^2\int_0^{+\infty}\frac{\rho^{n+1}}{\left(1+\rho^2\right)^n}d\rho$$
et 
$$J(\lambda)=\lambda^{n+2}(n-2)^2\sigma_{n-1} \int_0^{+\infty} \frac{\tau^{n-3}}{\left(\lambda^2+\tau^2\right)^{n}} v(\tau)d\tau\,.$$

Remarquons qu'il est facile de justifier nos calculs en les effectuant d'abord pour la fonction de Green conforme avec pôle en $x$
et condition  de Dirichlet sur le bord d'une grande boule géodésique $B(x,R)$ et en faisant ensuite tendre $R$ vers l'infini. De plus, le fait que $g=\bar g$ au dehors d'un voisinage $\mcO$ compact de $x$ et que
$$\int_{M\setminus \mcO} \scal_g G_x^2<+\infty$$ implique que
$$v(\tau)=O(r^{n-2}) \,\,\mathrm{et} \,\, \int_1^{+\infty}\frac{v(\tau)}{\tau^{n-1}}d\tau<+\infty,$$
et donc les quantités ci dessus sont bien finies.
\subsubsection{Asymptotiques}

Pour pouvoir faire un développement asymptotique des quantités $J(\lambda)$ et $ \|f_\lambda\|_{\frac{2n}{n-2}}^{\frac{2n}{n-2}}=:I(\lambda)$ lorsque $\lambda$ tend vers $0$, nous devons utiliser l'allure du noyau de Green près de $x$.
Celui-ci se trouve dans l'article de Lee-Parker (\cite[lem. 6.4]{LP}
\begin{prop}En utilisant les coordonnées normales $v=(v_1,\hdots , v_n)$ et $s=\sqrt{\sum_j v_j^2}$ de $\bar g$, on a :
\begin{enumerate}[i)]
\item Si $n\in \{3,4,5\}$ ou si $g$ est localement conformément plate au voisinage de $x$ alors
$$G_x(y)=\alpha_n\left(\frac{1}{s^{n-2}}+m+O(s)\right)$$
\item Si $n=6$ alors :
$$G_x(y)=\alpha_6\left(\frac{1}{s^{4}}-\frac{1}{288}|W|^2(x)\log s+O(1)\right)$$
\item Si $n\ge 7$ alors
$$G_x(y)=\frac{\alpha_n}{s^{n-2}}\left[1+\frac{s^2}{12(n-4)}\left(\frac{s^2}{12(n-6)} |W|^2(x)-\mathrm{Hess}_{\bar g}\scal_{\bar g}(v,v)\right)+O\left(s^5\right)\right]$$
\end{enumerate} 
\end{prop}
\`A partir de cette proposition, nous pouvons en déduire un développement asymptotique des fonctions $I$ et $J$.

Le plus simple est la fonction $J$ en tenant compte du fait que
$$\int_{\{s=R\}} \scal_{\bar g} d\sigma= -\frac{R^2}{2n}\sigma_{n-1} \left(\Delta_{\bar g}  \scal_{\bar g}\right)(x)+O(R^3)$$
on obtient facilement :
\begin{equation}\label{Jdas}
J(\lambda)=-\lambda^4 \frac{(n-2)^2\sigma_{n-1}}{192n(n-1)}|W|^2(x)\int_0^{+\infty} \frac{\rho^{n+1}}{\left(1+\rho^2\right)^n}d\rho+O\left(\lambda^5\right).
\end{equation}

De plus si la métrique $g$ est localement conformément plate au voisinage de $x$, alors on peut supposer que la métrique $\bar g$ est plate au voisinage de $x$ et dans ce cas on a
\begin{equation}\label{JdasW}
J(\lambda)=O\left(\lambda^{n+2}\right).
\end{equation}

Concernant la fonction $I$ on obtient :
\begin{enumerate}[i)]
\item Si $n\in \{3,4,5\}$ ou si $g$ est localement conformément plate au voisinage de $x$ alors
$$I(\lambda)=\omega_n A_{n,2}^{\frac{n}{n-2}}+m\lambda^{n-2}\frac{n}{n-2}\omega_n\int_0^\infty\frac{2n\rho^{2n-1}}{\left(1+\rho^2\right)^{n+1}}d\rho +o\left(\lambda^{n-2}\right)$$
\item Si $n=6$ alors on obtient
$$I(\lambda)=\omega_6 A_{6,2}^{\frac{3}{2}}-\frac{|W|^2(x)}{192}\lambda^4\log(\lambda)\int_0^\infty \frac{12\omega_6\rho^{11}}{\left(1+\rho^2\right)^{7}}d\rho+O\left(\lambda^4\right)$$
\item Si $n\ge 7$ alors
\begin{equation}\label{iii}I(\lambda)=\omega_n A_{n,2}^{\frac{n}{n-2}}+\frac{|W|^2(x)}{48(n-2)(n-6)}\lambda^4\int_0^\infty \frac{2n\omega_n\rho^{n+5}}{\left(1+\rho^2\right)^{n+1}}d\rho+O\left(\lambda^5\right).\end{equation}
\end{enumerate} 
Où 
\begin{equation}\label{A}\omega_n A_{n,2}^{\frac{n}{n-2}}=\sigma_{n-1}\int_0^{\infty}\frac{\rho^{n-1}}{(1+\rho^2)^n}d\rho=2^{-n}\sigma_n\end{equation}

\subsection{Discussion finale} 
De ces développements asymptotiques, on en déduit que si la dimension de $M$ égale à $6$ alors
$W(x)=0$. Ceci étant vrai pour tout $x$, on en déduit que la métrique $g$ est localement conformément plate.

Lorsque $n\ge 7$, nous devons comparer 
$$I(\lambda)^{1-\frac2n}$$ et
$$K(n,2)^2\left(\omega_n B_{n,2}-J(\lambda)\right).$$
Puisque
$$I(\lambda)=A_n+c_n|W|^2(x) \lambda^4+O\left(\lambda^5\right)$$
où les constantes $A_n$ et $c_n$ sont tirées de (\ref{iii}) on obtient
$$I(\lambda)^{1-\frac2n}=A_n^{1-2/n}+\left(1-\frac2n\right)c_nA_n^{-2/n} |W|^2(x)\lambda^4 +O\left(\lambda^5\right).$$
De plus on sait que 
$$K(n,2)^2=\sigma_n^{-2/n}\frac{4}{n(n-2)}\,\mathrm{et}\,\, \mathrm{que}\,\, A_n^{1-2/n}=K(n,2)^2\omega_n B_{n,2}.$$
De la même façon, on écrit
$$-J(\lambda)=b_n  |W|^2(x)\lambda^4+O\left(\lambda^5\right),$$ où la constante $b_n$ est issue de (\ref{Jdas}).
Compte-tenu de l'identité 
$$\int_{0}^{+\infty} \frac{\rho^a}{(1+\rho^2)^b} d\rho=\frac12\frac{\Gamma\left(\frac{a+1}{2}\right)\,\Gamma\left(n-\frac{a+1}{2}\right)}{\Gamma\left(b\right)},$$
on calcule
\begin{equation}\label{LHS}\begin{split}
\left(1-\frac2n\right)c_nA_n^{-2/n}&=\frac{n-2}{n}4 \sigma_n^{-2/n}\sigma_{n-1}\frac{1}{48\,(n-2)\,(n-6)}\frac{\Gamma\left(\frac{n}{2}+3\right)\Gamma\left(\frac{n}{2}-2\right)}{\Gamma\left(n+1\right)}\\
&= \sigma_n^{-2/n}\sigma_{n-1}\frac{1}{12\,n\,(n-6)}\frac{\left(\frac n2+2\right)\,\left(\frac n2+1\right)}{n\left(\frac{n}{2}-2\right)}\frac{\Gamma\left(\frac{n}{2}+1\right)\Gamma\left(\frac{n}{2}-1\right)}{\Gamma\left(n\right)}\\
&=\sigma_n^{-2/n}\sigma_{n-1}\frac{\Gamma\left(\frac{n}{2}+1\right)\Gamma\left(\frac{n}{2}-1\right)}{2\Gamma\left(n\right)}\frac{(n+4)(n+2)}{12n^2(n-4)(n-6)}
\end{split}
\end{equation}
et
\begin{equation}\label{RHS}\begin{split}
K(n,2)^2 b_n&= \sigma_n^{-2/n}\frac{4}{n(n-2)}\frac{\sigma_{n-1}(n-2)^2}{192 n(n-1)}\frac{\Gamma\left(\frac{n}{2}+1\right)\Gamma\left(\frac{n}{2}-1\right)}{2\Gamma\left(n\right)}\\
&=\sigma_n^{-2/n}\sigma_{n-1}\frac{\Gamma\left(\frac{n}{2}+1\right)\Gamma\left(\frac{n}{2}-1\right)}{2\Gamma\left(n\right)}\frac{n-2}{48n^2(n-1)}\,.
\end{split}
\end{equation}
Les calculs \ref{LHS} et \ref{RHS} montrent que (si $n>6$)
$$\left(1-\frac2n\right)c_nA_n^{-2/n}>K(n,2)^2 b_n \,.$$
Et on en déduit que la variété est localement conformément plate.

On en déduit que   $n\in \{3,4,5\}$ ou   que $g$ est localement conformément plate. Dans tous les cas (avec \ref{Jdas} ou \ref{JdasW}), on obtient alors que le terme de masse $m$ apparaissant dans le développement asymptotique du noyau de Green est forcément négatif ou nul.

Lorsque $n\ge 7$ et que  $g$ est localement conformément plate, le fait que la courbure scalaire de $g$ soit positive ou nulle, implique d'après \cite{SY} que l'on a conclusion du théorème. En effet, Schoen et Yau montre d'abord que dans ces dimensions, l'application développement de la structure localement conformément plate réalise un difféomorphisme conforme du  revêtement universel de $M$  vers un ouvert de la sphère \cite[Prop. 3.3 iii)]{SY}. D'après la seconde partie de l'argumentation de Schoen et Yau ceci implique que la masse est positive ou nulle en tout point et qu'elle est nulle si et seulement si $M$ est égale à son revêtement universel \cite[Prop. 4.4]{SY}.

Les autres cas découlent du théorème de masse positive généralisée et de l'adaptation des raisonnements de E.Witten  et de P. Jammes à notre cadre (\cite{W}, \cite{PT}, \cite{jammes}). 

Nous expliquons comment démontrer que la masse en $x$ est positive ou nulle lorsque $M$ est une variété spin que $n\in\{3,4,5\}$ ou que $(M^n,g)$ est localement conformément plate.
 Une adaptation de l'argument présenté ici et l'argument de P. Jammes permettra de démontrer que si $n$ est pair et si $g$ est localement conformément plat alors la masse est positive ou nulle.

On fixe donc $x$ et pour $\ve>0$ on considère la métrique
$$\widehat g_\ve =\left(\alpha_n^{-1}G_x+\ve\right)^{\frac{4}{n-2}}g$$
c'est une métrique complète sur $M\setminus\{x\}$ à courbure scalaire positive ou nulle et avec un bout asymptotiquement euclidien. Cette métrique étant conforme à $g$, elle vérifie aussi l'inégalité de Yamabe-Sobolev (\ref{YS}).
De plus la métrique $\widehat g_\ve$ admet dans des coordonnées stéréographiques autour de $x$ un développement limité :
$$\widehat g_\ve =\left(1+\frac{m+\ve}{t^{n-2}}\right)\eucl+O\left(t^{1-n}\right).$$
On note $\Sigma$ le fibré des spineurs sur $(M\setminus\{x\}, \widehat g_\ve )$
et on introduit l'espace de Hilbert
$H^1_0(M\setminus\{x\},\Sigma)$ obtenu en complétant 
l'espace $C^\infty_0(M\setminus\{x\},\Sigma)$ pour la norme associée à 
$$\sigma\mapsto \int_{M\setminus\{x\}} |\slashed D\sigma|^2 d\!\vol_{ \widehat g_\ve}=\int_{M\setminus\{x\}}\left[ |\nabla\sigma|^2+\frac14 \scal_{ \widehat g_\ve}|\sigma|^2\right]d\!\vol_{ \widehat g_\ve}$$
où on travaille avec la métrique $ \widehat g_\ve$.

Il est facile de voir que l'inégalité de Sobolev (\ref{YS}) implique que
$$\sigma \in H^1_0(M\setminus\{x\},\Sigma)\Rightarrow \nabla\sigma\in L^2\,\, \mathrm{et}\,\, \sigma\in L^{\frac{2n}{n-2}}$$ 
 De plus un argument facile basé sur des fonctions de coupure bien choisi permet de démontrer que si 
 $$\slashed D\sigma\in L^2\,\, \mathrm{et}\,\, \sigma\in L^2\Rightarrow \sigma \in H^1_0(M\setminus\{x\},\Sigma).$$
 
 Si $\sigma$ est un spineurs parallèle de norme $1$ sur l'espace euclidien, on en déduit par transplantation un spineur $\tilde\sigma$ sur $M\setminus\{x\}$ à support dans un voisinage de $x$ qui est asymptotiquement parallèle.
En cherchant le minima de la fonctionnelle définie sur $H^1_0(M\setminus\{x\},\Sigma)$ par
$$\xi\mapsto \frac12\int_{M\setminus\{x\}} |\slashed D\xi|^2 d\!\vol_{ \widehat g_\ve}-
\mathrm{Re}\left( \int_{M\setminus\{x\}} \langle \xi, \slashed D\tilde\sigma\rangle d\!\vol_{ \widehat g_\ve}\right),$$ 
on trouve $\xi\in H^1_0(M\setminus\{x\},\Sigma)$ vérifiant :
$$\slashed D^2\xi=\slashed D\tilde \sigma.$$
On en déduit que $\slashed D\xi \in H^1_0(M\setminus\{x\},\Sigma)$ ainsi le spineur
défini par
$$\widehat \sigma=\tilde \sigma-\slashed D\xi$$
est harmonique  $$\slashed D\widehat \sigma=0$$ et de plus 
$\nabla \widehat\sigma\in L^2$, car $\nabla \tilde \sigma\in L^2$ et $\slashed D(\slashed D\xi)$ et $\slashed D\xi$ sont dans $L^2$ donc $\nabla\slashed D\xi\in L^2$.

Le fait que $(M\setminus\{x\},\widehat g_\ve)$ soit complet permet de justifier l'intégration par parties de Witten et d'obtenir
$$\frac{\sigma_{n-1}}{4}(m+\ve)=\int_{M\setminus\{x\}}\left[ |\nabla\widehat\sigma|^2+\frac14 \scal_{ \widehat g_\ve}|\widehat\sigma|^2\right]d\!\vol_{ \widehat g_\ve}.$$

Ainsi pour tout $\ve > 0$, on a $m+\ve\ge 0$ donc la masse est positive ou nulle.
Dans le cadre du théorème C, on sait alors que la masse est nulle.
 Alors les métriques $ \widehat g_\ve$ et $\widehat g_0 =\left(\alpha_n^{-1}G_x\right)^{\frac{4}{n-2}}g$ étant conformément équivalente, en utilisant l'invariance conforme de l'opérateur de Dirac et en passant à la limite $\ve\to 0$ on obtient :
$$0=\int_{M\setminus\{x\}} |\nabla\widehat\sigma|^2 d\!\vol_{ \widehat g_0}.$$
Où toutes les quantités sont calculées par rapport à la métrique non complète $\widehat g_0$.
On en déduit que le spineur $\widehat\sigma$ est parallèle. On en déduit que le fibré des spineurs de $(M\setminus\{x\},\widehat g_0)$ est trivialisé par une base de spineurs parallèle. On en déduit donc que la métrique
$(M\setminus\{x\},\widehat g_0)$ est plate sans holonomie. On peut alors conclure avec les arguments de Schoen et Yau.

 \end{document}